\documentclass[8pt,bezier]{article}
\usepackage{amsmath,amssymb,amsfonts,euscript}
\usepackage{algorithmic, algorithm}
\usepackage{graphicx}

\usepackage[ansinew]{inputenc}
\usepackage{graphicx}
\usepackage{color}
\usepackage{mathrsfs}

\usepackage[colorlinks]{hyperref}
\usepackage[active,new,noold,marker]{xrcs}
\catcode`\=13
\def{$\bowtie$}

\usepackage{eurosym}
\usepackage{lscape} 

\newtheorem{prethm}{{\bf Theorem}}

\newenvironment{thm}{\begin{prethm}{\hspace{-0.5
               em}{\bf}}}{\end{prethm}}

\newtheorem{prepro}{{\bf Theorem}}

\newtheorem{preprop}{{\bf Proposition}}

\newtheorem{precor}{{\bf Corollary}}

\newtheorem{preconj}{{\bf Conjecture}}

\newenvironment{conj}{\begin{preconj}{\hspace{-0.5
               em}{\bf}}}{\end{preconj}}

\newtheorem{predefi}{{\bf Definition}}

\newtheorem{preremark}{{\bf Remark}}

\newenvironment{remark}{\begin{preremark}\rm{\hspace{-0.5
               em}{\bf}}}{\end{preremark}}

\newtheorem{preexample}{{\bf Example}}

\newtheorem{prelem}{{\bf Lemma}}

\newtheorem{prelam}{{\bf Theorem}}

\newenvironment{lam}{\begin{prelam}{\hspace{-0.5
               em}{\bf}}}{\end{prelam}}

\newtheorem{preprob}{{\bf Problem}}

\newenvironment{prob}{\begin{preprob}{\hspace{-0.5
               em}{\bf.}}}{\end{preprob}}

\newtheorem{preproof}{{\bf Proof}}

\newtheorem{preali}{{\bf Proof of Theorem 1.}}

\newtheorem{prealii}{{\bf Proof of Theorem 2.}}

\newtheorem{prealiii}{{\bf Proof of Theorem 3.}}

\title{ Algorithmic Complexity of Weakly Semiregular Partitioning and the Representation Number}

\author{\bf\small Arash Ahadi  ${}^{a}$, {Ali Dehghan}
${}^{b}$ \thanks{Corresponding author. \newline
${}^{a}$ Department of Mathematical Sciences, Sharif University of Technology, Tehran, Iran
\newline ${}^{b}$ Systems and Computer Engineering Department, Carleton University, Ottawa,   Canada
\newline ${}^{c}$ Department of Mathematics, University of Western Ontario, London, Ontario, Canada
\newline E-mail Addresses:
arash$_{-}$ahadi@mehr.sharif.edu (Arash Ahadi)
alidehghan@sce.carleton.ca (Ali Dehghan),
mmollaha@uwo.ca (Mohsen Mollahajiaghaei) .} , {Mohsen Mollahajiaghaei}
${}^{c}$ ~}

\date{}
\begin{document}
\maketitle

\begin{abstract}
{\small \noindent
A graph $G$ is {\it weakly semiregular} if there are two numbers $a,b$, such that the degree of every vertex is $a$ or $b$.
The {\it weakly semiregular number} of a graph $G$, denoted by $wr(G)$, is the minimum number of subsets into which the edge set of $G$ can be partitioned
so that the subgraph induced by each subset is a weakly semiregular graph.
We present a polynomial time algorithm to determine whether the weakly semiregular number of a given tree is two.
On the other hand, we show that determining whether $ wr(G) = 2 $ for a given bipartite graph $ G $ with at most three numbers in its degree set is {\bf NP}-complete.
Among other results, for every tree $T$, we show that $wr(T)\leq 2\log_2 \Delta(T) + \mathcal{O}(1)$, where $\Delta(T)$ denotes the  maximum
degree of $T$.\\
A graph $G$ is a {\it $[d, d +s]$-graph} if the degree of every vertex of $G$ lies in the interval $[d, d +s]$. A $[d, d +1]$-graph is said to be {\it semiregular}. The {\it semiregular number} of a graph $G$, denoted by $sr(G)$, is the minimum number of subsets into which the edge set of $G$ can be partitioned
so that the subgraph induced by each subset is a semiregular graph. We prove that the semiregular number of a  tree $T$ is  $ \lceil \frac{\Delta(T)}{2}\rceil$.
On the other hand, we show that determining whether $ sr(G) = 2 $ for a given bipartite graph $ G $ with $\Delta(G)\leq 6$  is {\bf NP}-complete.\\
In the second part of the work, we consider the representation number. A graph $G$ has a {\it representation modulo $r$} if there exists an injective map $\ell: V (G) \rightarrow
\mathbb{Z}_r$ such that vertices $v$ and $u$ are adjacent if and only if $|\ell(u) -\ell(v)|$ is
relatively prime to $r$.
The {\it representation number}, denoted by $rep(G)$, is the smallest $r$ such that $G$ has a representation
modulo $r$. Narayan and Urick conjectured that the determination of $rep (G)$ for an arbitrary graph $G$ is a  difficult problem \cite{narayan2007representations}. In this work, we confirm this conjecture and show that
if $\mathbf{NP\neq P}$, then for any $\epsilon >0$, there is no polynomial time
$(1-\epsilon)\frac{n}{2}$-approximation algorithm for the computation of
representation number of regular graphs with $n$ vertices.
}

\begin{flushleft}

\end{flushleft}

\end{abstract}

\section{Introduction}
\label{}

The paper consists of two parts. In the first part,
we consider the problem of partitioning the edges of
a graph into  regular and/or locally irregular subgraphs. In this part, we present some polynomial time algorithms and {\bf NP}-hardness results. In the second part of the work, we focus on the representation number of graphs. It was conjectured that the determination of $rep (G)$ for an arbitrary graph $G$ is a  difficult problem \cite{narayan2007representations}. In this part, we confirm this conjecture and show that
if $\mathbf{NP\neq P}$, then for any $\epsilon >0$, there is no polynomial time
$(1-\epsilon)\frac{n}{2}$-approximation algorithm for the computation of
representation number of regular graphs with $n$ vertices.

\section{Partitioning the edges of graphs}

In 1981, Holyer in \cite{MR635429} focused on the computational complexity of edge partitioning problems and
 proved that for each $t$, $t \geq 3$, it is {\bf NP}-complete to decide whether a given graph
can be edge-partitioned into subgraphs isomorphic to the complete graph $K_t$.
Afterwards, the complexity of edge partitioning problems have been studied extensively by several
authors, for instance see \cite{dehghan2016colorful, diwan2015complexity,  MR1460720}.
Nowadays, the computational complexity of edge partitioning problems is a well-studied area of graph theory and computer science.
For more information we refer the reader to a survey
on graph factors and factorization by Plummer \cite{plummer2007graph}.

If we consider the Holyer problem for a family $\mathcal{G}$ of graphs instead of a fixed graph
then, we can discover interesting problems.
For a family $\mathcal{G}$ of graphs, a  $\mathcal{G}$-decomposition of a graph $G$ is a partition of the edge set of $G$ into
subgraphs isomorphic to members of $\mathcal{G}$.
Problems of $\mathcal{G}$-decomposition of graphs have received a considerable attention,
for example, Holyer proved that it is {\bf NP}-hard to edge-partition a graph into the minimum number of complete subgraphs \cite{MR635429}. To see more examples of  $\mathcal{G}$-decomposition of graphs see \cite{bensmail2014partitions, chung1981decomposition, kirkpatrick1983complexity}.

\subsection{Related works and motivations}

We say that a graph is  {\em locally irregular} if its adjacent vertices have
distinct degrees and a graph is {\em regular} if each vertex of the graph has the same degree.
In 2001, Kulli {\it et al.} introduced an interesting parameter for the partitioning of the edges of a graph \cite{kulli}.
The {\it regular number} of a graph $G$, denoted by $reg(G)$, is the minimum number of subsets into which the edge set of $G$ can be partitioned so that the subgraph induced by each subset is regular.
The  {\it edge chromatic number} of a graph, denoted by $\chi '(G)$, is the minimum size of a partition of the edge set into $1$-regular subgraphs and also, by Vizing's theorem \cite{MR0180505}, the edge chromatic number of a graph $G$ is equal to either $ \Delta(G) $ or $ \Delta(G) +1 $, therefore the regular number problem is a generalization for the  edge chromatic number and we have the following bound: $ reg(G)\leq \chi '(G) \leq \Delta(G) +1$.
It
was asked  in \cite{kulli2} to determine whether  $reg(G)\leq \Delta(G)$ holds for all
connected graphs.

\begin{conj} \label{III0} [ The degree bound {\it \cite{kulli2}}]
For any connected graph $G$, $reg(G)\leq \Delta(G)$.
\end{conj}

It was shown in \cite{dehghan2014complexity} that not only  there exists a counterexample for the above-mentioned bound but also for a given connected graph $G$ decide whether $reg(G)=\Delta(G)+1$ is
{\bf NP}-complete.
Also, it was shown that the computation of the regular number for a given connected bipartite graph $G$ is {\bf NP}-hard \cite{dehghan2014complexity}. Furthermore, it was proved that determining whether $ reg(G) = 2 $ for a given connected $3$-colorable graph $ G $  is {\bf NP}-complete \cite{dehghan2014complexity}.

On the other hand, Baudon {\it et al.} introduced the notion
of edge partitioning into locally irregular subgraphs\cite{baudon2013decomposing}.
In this case, we want to partition the edges of the graph $G$ into
locally irregular subgraphs, where by a partitioning of the graph $G$ into $k$ locally
irregular subgraphs
 we refer to a partition $E_1, \ldots, E_k$ of $E(G)$ such
that the graph  $G[E_i]$ is locally irregular for every $i$, $i =1,\ldots,k$.
The {\it irregular chromatic index} of $G$, denoted by $ \chi'_{irr}$,
is the minimum number $k$ such that the graph $G$ can be partitioned into $k$ locally
irregular subgraphs. Baudon {\it et al.}
characterized all graphs which cannot be partitioned into locally
irregular subgraphs and  call them exceptions \cite{baudon2013decomposing}.
Motivated by the 1-2-3-Conjecture, they conjectured that apart from these exceptions all other connected
graphs can be partitioned into three locally irregular subgraphs \cite{baudon2013decomposing}. For more information about the 1-2-3-Conjecture and its variations, we refer the reader to a survey on the 1–2–3 Conjecture and related
problems by Seamone \cite{survey} (see also \cite{ MR3529282, ahadi, Barme2016, MR3589724, MR3512668,  MR2654258, MR3022926}).

\begin{conj} \label{III1} {\it \cite{baudon2013decomposing}}
For every non-exception graph $G$, we have $\chi'_{irr}(G)\leq 3$.
\end{conj}

Regarding the above-mentioned conjecture, Bensmail {\it et al.} in \cite{bensmail2016decomposing} proved that every bipartite graph $G$ which is not an odd
length path satisfies $\chi'_{irr}(G)\leq 10$. Also, they proved that if $G$ admits a partitioning
into locally irregular subgraphs, then $\chi'_{irr}(G)\leq 328 $. Recently, Lu{\v{z}}ar {\it et al.} improved the previous bound for bipartite graphs and general graphs to 7 and 220, respectively \cite{luvzar2016new}. For more information about this conjecture see \cite{MR3512653}.

Regarding the complexity of edge partitioning into locally irregular subgraphs,  Baudon {\it et al.} in  \cite{bensmail2013complexity} proved that the problem of determining
the irregular chromatic index of a graph can be handled in linear
time when restricted to trees. Furthermore, in \cite{bensmail2013complexity}, Baudon {\it et al.} proved that determining whether a given planar graph $G$
can be partitioned into two locally irregular subgraphs is {\bf NP}-complete.

In 2015, Bensmail and Stevens considered the problem of partitioning the edges of graph into
some subgraphs, such that in each subgraph every component is either
regular or locally irregular \cite{bensmail2014edge}. The {\it regular-irregular chromatic index} of graph $G$,
 denoted by $ \chi'_{reg-irr}$,
is the minimum number $k$ such that $G$ can be partitioned into $k$ subgraphs, such that each component  of every
subgraph is locally
irregular or regular \cite{bensmail2014edge}.
They conjectured that the edges of every graph can be
partitioned into at most two subgraphs, such that each component of every subgraph
 is regular or locally irregular \cite{ bensmail2014decomposing, bensmail2014edge}.

\begin{conj} \label{III2}  {\it \cite{bensmail2014decomposing, bensmail2014edge}}
For every  graph $G$, we have $\chi'_{reg-irr}(G)\leq 2$.
\end{conj}

Recently, motivated by Conjecture \ref{III1} and Conjecture \ref{III2}, Ahadi {\it et al.} in \cite{Ahadi-deh} presented  the following conjecture.
With Conjecture \ref{III3} they weaken Conjecture \ref{III1} and strengthen Conjecture \ref{III2}.

\begin{conj} \label{III3}  {\it \cite{Ahadi-deh}}
Every graph can be partitioned into 3 subgraphs, such that each subgraph
is locally irregular or regular.
\end{conj}

Note that in Conjecture \ref{III3},  each subgraph (instead of each component of every subgraph) should be locally irregular or regular.
Also, note that it was shown that
deciding whether a given planar bipartite graph $G$ with maximum degree three can be partitioned into at most
two subgraphs such that each subgraph is regular or locally irregular is {\bf NP}-complete \cite{Ahadi-deh}.

In \cite{Ahadi-deh}, Ahadi {\it et al.} considered the problem of partitioning the edges into locally regular subgraphs.
We say that a graph $G$ is {\it locally regular} if each component of $G$ is regular (note that a regular graph is locally regular but the converse does not hold).
The {\it  regular chromatic index} of a graph $G$ denoted by $\chi'_{reg}$ is the minimum number of subsets into which the
edge set of $G$ can be partitioned so that the subgraph induced by each subset is locally regular.
From the definitions of locally regular and regular graphs we have the following bound: $\chi'_{reg}(G) \leq reg(G)\leq  \Delta(G)+1$. It was shown that every graph $G$ can be partitioned into $\Delta(G)$
subgraphs such that each subgraph is locally regular and this bound is sharp for trees \cite{Ahadi-deh}.

\begin{lam} \cite{Ahadi-deh}
Every graph $G$ can be partitioned into $\Delta(G)$
subgraphs such that each subgraph is locally regular and this bound is sharp for trees.
\end{lam}

In conclusion, we can say that the problem of partitioning the edges of graph into regular and/or locally irregular subgraphs is an active area in graph theory and computer science.
What can we say about the edge decomposition problem if
 we require that each subgraph (instead of each component of every subgraph)
should be a graph with at most $k$  numbers in its degree set.
With this motivation in mind, we investigate
the problem of partitioning the edges of graphs into subgraphs such that each subgraph has at most two numbers in its degree set. In this work, we consider partitioning into weakly semiregular and semiregular subgraphs.

\subsection{Weakly semiregular graphs}

A graph $G$ is {\it weakly semiregular} if there are two numbers $a,b$, such that the degree of every vertex is $a$ or $b$.
The {\it weakly semiregular number} of a graph $G$, denoted by $wr(G)$, is the minimum number of subsets into which the edge set of $G$ can be partitioned so that the subgraph induced by each subset is weakly semiregular.
This parameter is well-defined for any graph $G$ since one can
always partition the edges into $1$-regular subgraphs.
Throughout the paper, we say that a graph $G$ is $(a,b)$-graph if the degree of every vertex is $a$ or $b$ (in other words, if the degree set of the graph $G$ is $\{a,b\}$).

\begin{remark}
There are infinitely many values of $\Delta$ for which the graph $G$ might be chosen so that $wr(G)\geq \log_3 \Delta(G)$.
Assume that $G$ is a graph such that for each $i$, $1\leq i \leq \Delta$, there is a vertex with degree $i$ in that graph.
Also, let $E_1,E_2,\ldots,E_{wr(G)}$ be a weakly semiregular partitioning for the edges of that graph.
The degree set of the subgraph $G_i = (V,E_i)$ has at most three elements. By adding $wr(G)$
such degree sets, one corresponding to each subset $E_i$, we get a degree set
that contains at most $3^{wr(G)}$ elements. Hence, the degree set of the graph $G$ contains
at most $3^{wr(G)}$ elements. This completes the proof.
\end{remark}

In this work, we focus on the algorithmic aspects of weakly semiregular number.
We present a polynomial time algorithm to determine whether the weakly semiregular number of a given tree is two.

\begin{thm}\label{T1}\\
(i) There is an $\mathcal{O}(n^2)$ time algorithm  to determine whether the weakly semiregular number of a given tree is two, where $n$ is the number of vertices in the tree.\\
(ii) Let $c$ be a constant, there is a polynomial  time algorithm  to determine whether the weakly semiregular number of a given tree is at most $c$.\\
(iii) For every tree $T$, $wr(T)\leq 2\log_2 \Delta(T) + \mathcal{O}(1)$.
\end{thm}

\begin{remark}\label{R1}
If $G$ is a graph with $\Delta(G)\leq 4$, then $wr(G)\leq 2$. If the graph $G$ is not regular, then consider two copies of the graph $G$ and
for each vertex $v$ with degree less than 4, join the vertex $v$ in the first copy of $G$ to the vertex $v$ in the second copy of the graph $G$. By repeating this procedure we can obtain a 4-regular graph $G'$.
A subgraph $F$ of a graph $H$ is called a factor of $H$ if $F$ is a spanning
subgraph of $H$. If a factor $F$ has all of its degrees equal to $k$, it is called a $k$-factor.
A $k$-factorization for a graph $H$ is a partition of the edges
into disjoint $k$-factors. For $k\geq 1$, every $2k$-regular graph admits a $2$-factorization \cite{petersen1891theorie},
thus the graph $G'$ can be partitioned into two 2-regular graphs $G_1'$ and $G'_2$. Let $f:E(G')\rightarrow\{1,2\}$ be a function such that $f(e)=1$ if and only if $e\in E(G_1')$. One can see that the function $f$ can partition the edges of  the graph $G$ into two (1,2)-graphs. Therefore, $wr(G)\leq 2$. This completes the proof.
\end{remark}

If $G$ is a graph with at most two numbers in its degree set, then its weakly semiregular number is one.
On the other hand, if $\Delta\leq 4$ by Remark \ref{R1}, the weakly semiregular number of the graph is at most two.
We show that determining whether $ wr(G) = 2 $ for a given bipartite graph $ G $ with $\Delta(G)=6$ and at most three numbers in its degree set,  is {\bf NP}-complete.\\

\begin{thm}\label{T2}
Determining whether $ wr(G) = 2 $ for a given bipartite graph $ G $ with $\Delta(G)=6$ and at most three  numbers in its degree set,  is {\bf NP}-complete.
\end{thm}

\subsection{Semiregular graphs}

A graph $G$ is a {\it $[d, d +s]$-graph} if the degree of every vertex of $G$ lies in the interval $[d, d +s]$. A $[d, d +1]$-graph is said to be {\it semiregular}. Semiregular graphs are an important family of graphs and their properties have been studied extensively, see for instance \cite{balbuena2010connectivity, balbuena2010some}.
The {\it semiregular number} of a graph $G$, denoted by $sr(G)$, is the minimum number of subsets into which the edge set of $G$ can be partitioned
so that the subgraph induced by each subset is semiregular. We prove that the semiregular number of a tree $T$ is $ \lceil \frac{\Delta(T)}{2}\rceil$.
On the other hand if $\Delta\leq 4$ by Remark \ref{R1}, the semiregular number of a graph is at most two.
We show that determining whether $ sr(G) = 2 $ for a given bipartite graph $ G $ with $\Delta(G)\leq 6$,  is {\bf NP}-complete.

\begin{thm}\label{T3}\\
(i) Let $T$ be a tree, then $sr(T) = \lceil\frac{\Delta(T)}{2}\rceil$.\\
(ii) Let $G$ be a graph, then $sr(G) \leq \lceil\frac{\Delta(G)+1}{2}\rceil$.\\
(iii) Determining whether $ sr(G) = 2 $ for a given bipartite graph $ G $ with $\Delta(G)\leq 6$,  is {\bf NP}-complete.
\end{thm}

Every semiregular graph is a weakly semiregular graph, thus by the above-mentioned theorem, we have the following bound:

\begin{equation}
wr(G) \leq sr(G) \leq \lceil\frac{\Delta(G)+1}{2}\rceil
\end{equation}

\subsection{Partitioning into locally irregular and weakly semiregular subgraphs}

In \cite{bensmail2014edge}
Bensmail  and   Stevens considered the outcomes on Conjecture \ref{III1} of
allowing components isomorphic to the complete graph $K_2$, or more generally regular components.
In fact their investigations are  motivated by the following question: "How easier can
 Conjecture \ref{III1} be tackled if we allow a locally irregular partitioning to induce connected
components isomorphic to the complete graph $K_2$?"
They conjectured that the edges of every graph can be
partitioned into at most two subgraphs, such that each component of every subgraph
is regular or locally irregular \cite{ bensmail2014edge}.
Motivated by this conjecture we pose the following conjecture.
Note that in Conjecture \ref{III6},  each subgraph (instead of each component of every subgraph) should be locally irregular or weakly semiregular.

\begin{conj} \label{III6}
Every graph can be partitioned into 3 subgraphs, such that each subgraph
is locally irregular or weakly semiregular.
\end{conj}

Note that if  Conjecture \ref{III1} or Conjecture \ref{III3}  is true, then Conjecture \ref{III6} is true. Also, if every graph can be partitioned into 2 subgraphs such that each component of every subgraph is a locally irregular graph or $K_2$, then Conjecture \ref{III6} is true.
We conclude this section by the following hardness result.

\begin{thm}\label{T4}
Determining whether a given  graph $G$,
can be partitioned into 2 subgraphs, such that each subgraph
is locally irregular or weakly semiregular is {\bf NP}-complete.
\end{thm}

\subsection{Summary of results}

A summary of results and open problems on edge-partition problems are shown in Table 1.

\begin{table}[ht]
\small
\caption{Recent results on edge partitioning of graphs into subgraphs} 
\centering 
\begin{tabular}{l l  l l l} 
\hline\hline 
                 & $=2$ (for trees)  &  $=2$  & Upper bound  \\ [0.5ex] 
\hline 
Regular subgraphs &  {\bf P} (see \cite{kulli2}) & {\bf NP}-c  (see \cite{dehghan2014complexity})  &  $\Delta+1$ (see \cite{kulli2}) \\
Locally regular subgraphs &  {\bf P} (see \cite{Ahadi-deh}) & {\bf NP}-c  (see \cite{Ahadi-deh})  &  $\Delta $ (see \cite{Ahadi-deh}) \\
Weakly semiregular subgraphs &  {\bf P} (Th. \ref{T1}) & {\bf NP}-c (Th. \ref{T2}) &  $\lceil\frac{\Delta+1}{2}\rceil$ (Th. \ref{T3}) \\
Semiregular  subgraphs &  {\bf P} (Th. \ref{T3}) & {\bf NP}-c (Th. \ref{T3}) &  $\lceil\frac{\Delta+1}{2}\rceil$ (Th. \ref{T3}) \\
Locally irregular subgraphs &  {\bf P} (see \cite{bensmail2013complexity}) & {\bf NP}-c (see \cite{bensmail2013complexity}) &  3 (Conj.  \ref{III1}) \\

regular-irregular subgraphs &  Open  (see \cite{Ahadi-deh}) & {\bf NP}-c (see \cite{Ahadi-deh}) &  3 (Conj.  \ref{III3}) \\
regular-irregular components &  {\bf P} (see \cite{bensmail2014edge}) & P (Conj.  \ref{III2}) &  2 (Conj.  \ref{III2}) \\

\hline 
\end{tabular}
\label{table:nonlin} 
\end{table}

\section{Representation number}

A finite graph $G$ is said to be {\it representable modulo $r$}, if there exists an injective map $\ell: V (G) \rightarrow
\{0,1,\ldots,r-1\}$ such
that vertices $v$ and $u$ of the graph $G$ are adjacent if and only if $|\ell(u) -\ell(v)|$ is
relatively prime to $r$. The {\it representation number} of $G$,
denoted by $rep(G)$, is the smallest positive integer $r$ such that the graph $G$ has a representation
modulo $r$. In 1989, Erd\H{o}s and Evans introduced representation numbers and showed that every finite graph can be represented modulo some
positive integer \cite{erdos1989representations}. They
used representation numbers to give
a simpler proof of a result of Lindner et al. \cite{lindner1979orthogonal} that, any finite graph can be realized as an orthogonal Latin square graph (an orthogonal Latin square graph is one whose vertices can be labeled with Latin squares
of the same order and same symbols such that two vertices are adjacent if and only if the
corresponding Latin squares are orthogonal). The existence proof of Erd\H{o}s and Evans  gives an unnecessarily large upper bound for
the representation number \cite{erdos1989representations}.
During the recent years, representation numbers
have received considerable attention and
have been studied for various classes of graphs, see \cite{akhtar2012representation, MR1782042,  evans1994representations, MR3278261, narayan2003upper}.

Narayan and Urick conjectured that the determination of $rep (G)$ for an arbitrary graph $G$ is a  difficult problem \cite{narayan2007representations}.
In the  following theorem we discuss about the computational complexity
of $rep(G)$ for regular graphs.

\begin{thm} \label{T5}\\
(i) If $\mathbf{NP\neq P}$, then for any $\epsilon >0$, there is no polynomial time
$(1-\epsilon)\frac{n}{2}$-approximation algorithm for the
representation number of regular graphs with $n$ vertices.\\
(ii) For every $\epsilon > 0$ there is a   polynomial time $((1+\epsilon)\dfrac{en}{2})$-approximation algorithm for
computing $rep(G)$ where $G^c$ is a triangle-free  $r$-regular
graph.
\end{thm}

\section{Notation and tools}

All graphs considered in this paper are finite and undirected. If $G$ is a graph, then $V(G)$ and $E(G)$ denote the vertex set and the edge set of
$G$, respectively. Also, $\Delta(G)$ denotes the maximum degree of $G$ and simply denoted
by $\Delta$. For every $v\in V(G)$, $d_{G}(v)$ and $N_{G}(v)$ denote the degree of $v$ and the set of
neighbors of $v$, respectively. Also, $N[v]=N(v)\cup \{v\}$.
For a given graph $G$, we use $u\sim v$ if two vertices $u$ and $v$ are adjacent in
$G$.

The {\it degree sequence} of a graph
is the sequence of non-negative integers listing the degrees of the vertices of $G$. For
example, the complete bipartite graph $K_{1,3}$ has degree sequence $(1, 1, 1,  3)$, which contains two distinct elements:
1 and 3. The {\it degree set} $D$ of a graph $G$ is the set of distinct degrees of the vertices of $G$.
For $ k\in \mathbb{N} $, a {\it proper edge $k$-coloring} of $G$ is a function $c:
E(G)\longrightarrow \lbrace 1,\ldots , k \rbrace$, such that if $e,e'\in E(G)$ share a common endpoint,
then $c(e)$ and $c(e')$ are different.
The smallest integer $k$ such that
$G$ has a proper edge $k$-coloring is called the {\it edge chromatic number} of $G$ and denoted by $\chi '(G)$. By Vizing's theorem \cite{MR0180505}, the edge chromatic number of a graph $G$ is equal to either $ \Delta(G) $ or $ \Delta(G) +1 $. Those graphs
$G$ for which $\chi '(G)=\Delta(G)  $ are said to belong to $Class$ $1$, and the other to $Class$ $2$.

Let $G$ be a graph and $f$ be a non-negative integer-valued function on $V (G)$. Then a spanning
subgraph $H$ of $G$ is called an $f$-factor of $G$ if $d_H (v) = f (v)$, for all $v \in V (G)$.
Let $G$ be a graph  and let $f$, $g$ be mappings of $  V (G)$ into the non-negative integers.
A $(g, f )$-factor of $G$ is a spanning subgraph $F$ such that $g(v) \leq d_F(v) \leq f (v)$ for all $v\in V (G)$.
In 1985, Anstee  gave a polynomial time algorithm for the $(g, f )$-factor problem and his algorithm
either returns one of the factors in question or shows that none exists, in $\mathcal{O}(n^3)$ time \cite{Anstee}. Note that this complexity bound is independent of the functions $g$ and $f$. We will use form this algorithm in our proof.
We
follow \cite{MR1367739} for terminology and notation where they are not defined here.

\section{Proof of Theorem \ref{T1}}

(i) Let $T$ be an arbitrary tree. Any subgraph of a tree is a forest, so if $T$ can be partitioned into two weakly semiregular forests $T_1$ and $T_2$, then there are two numbers $\alpha , \beta$ (not necessary distinct) such that $T_1$ is a $(1,\alpha)$-forest and $T_2$ is a $(1,\beta)$-forest (note that a forest $T$ is a $(a, b)$-forest if the degree of every vertex
is $a$ or $b$). Without loss of generality, we can assume that $1 \leq \alpha \leq \beta \leq \Delta(T) \leq n$. Let $D$ be the degree set of $T$, we have $D \subseteq \{1,2,\alpha , \alpha+1 , \beta , \beta+1, \alpha+\beta\}$. So if $|D|\geq 8$, then the tree $T$ cannot be partitioned into two weakly semiregular forests.
On the other hand, one can see that if $|D|\leq 7$, then the number of possible cases for $(\alpha, \beta)$ is $\mathcal{O}(1)$.

In Algorithm 1, we present an $\mathcal{O}(n^2)$ time algorithm to check whether $T$ can be partitioned into
two weakly semiregular forests $T_1$ and $T_2$, such that the forest $T_1$ is $(1,\alpha)$-forest and the forest $T_2$ is $(1,\beta)$-forest.
If the algorithm returns NO, it means that $T$ cannot be partitioned and if it returns YES, it means that $T$ can be partitioned.

Here, let us to introduce some notation and state a few properties of Algorithm 1.
Suppose that $|V(T)| =n$ and choose an arbitrary vertex $v\in V(T)$ to be its root.
Perform a breadth-first search
algorithm from the vertex $v$. This defines a partition $L_{0}, L_{1}, \ldots , L_{h}$ of the vertices of $T$
where each part $L_i$ contains the vertices of $T$ which are at depth $i$ (at distance exactly $i$
from $v$). Let $p(x)$ denote the neighbor of the vertex $x$ on the $xv$-path, i.e. its parent. Also,
 let $\lbrace v_{1}=v,v_2,\ldots , v_{n} \rbrace  $ be a list of the vertices according to their distance from the root.
We use form this list of vertices in the algorithm. See Algorithm 1.

\begin{algorithm}
\algsetup{linenosize=\tiny}
\small
\caption{}
\begin{algorithmic}[1]

\STATE {{\bf Input:}} The tree $T$ and two numbers $\alpha, \beta$.
\STATE {{\bf Output:}}  Can $T$  be partitioned into
two weakly semiregular forests $T_1$ and $T_2$, such that $T_1$ is $(1,\alpha)$-forest and $T_2$ is $(1,\beta)$-forest.
\STATE Let $g:E(T) \rightarrow \{red, blue, free\}$ and put $g(e)\leftarrow free$ for all edges
\STATE Let $f:E(T) \rightarrow \{red, blue, free\}$ and put $f(e)\leftarrow free$ for all edges
\WHILE{ there is an edge $e$ such that $g(e)=free$}
  \STATE For any edge $e$, put $f(e)\leftarrow free$
  \STATE $s\leftarrow $YES
  \FOR{ $i=1$ to $i=n$}
    \IF {there is no labeling like $h$ for the set of edges $S_i=\{v_iv_j: j>i\}$ with the colors red and blue such that $|\{e: e\in S_i, h(e)=red\}\cup\{v_ip(v_i):   f(v_ip(v_i))=red, i>1\}|\in \{0 , 1 ,\alpha\}$, also $|\{e: e\in S_i, h(e)=blue\}\cup\{v_ip(v_i):  f(v_i)=blue, i>1\}|\in \{0 , 1 ,\beta\}$
    and for each edge $e\in S_i$, if $g(e)\neq free$, then $g(e)=h(e)$}
        \IF{ $g(v_ip(v_i))\neq free$}
        \RETURN NO
        \ENDIF
         \IF{ $f(v_ip(v_i))=blue$}
          \STATE $g(v_ip(v_i))\leftarrow red$
           \STATE $s\leftarrow $NO
          \STATE {\bf break} the for loop
        \ENDIF
        \IF{ $f(v_ip(v_i))=red$}
          \STATE $g(v_ip(v_i))\leftarrow blue$
          \STATE $s\leftarrow $NO
          \STATE {\bf break} the for loop
        \ENDIF
     \ENDIF
     \IF{$s=$YES}
    \STATE Label the set of edges $S_i=\{v_iv_j: j>i\}$ with the colors red and blue such that $|\{e: e\in S_i, h(e)=red\}\cup\{v_ip(v_i): \text{ if } f(v_i)=red\}|\in \{0 , 1 ,\alpha\}$, also $|\{e: e\in S_i, h(e)=blue\}\cup\{v_ip(v_i): \text{ if } f(v_i)=blue\}|\in \{0 , 1 ,\beta\}$
    and for each edge $e\in S_i$, if $g(e)\neq free $, then $g(e)=h(e)$
    \STATE For each edge $e$ in $\{v_iv_j: j>i\}$ assign the label of $e$ to the variable $f(e)$
    \ENDIF
  \ENDFOR
  \IF{$s=$YES}
   \RETURN YES
  \ENDIF
\ENDWHILE
\RETURN NO
\end{algorithmic}
\end{algorithm}

{\bf Sketch of Algorithm 1}\\
The algorithm starts from the vertex $v_1$ and labels the set of edges $S_1=\{v_1v_j: j>1\}$ with labels {\it red} and {\it blue} such that the number of edges in $S_1$ with label {\it red} is 0 or 1 or $\alpha$ and the number of edges in $S_1$ with label {\it blue} is 0 or 1 or $\beta$. The algorithm saves the partial labeling in $f$. Next, at step $i$, $i>1$ of the {\it for loop}, the algorithm labels the set of edges $S_i=\{v_iv_j: j>i\}$ with labels {\it red} and {\it blue} such that the number of edges in $S_i \cup \{v_ip(v_i) \}$ with label {\it red} is 0 or 1 or $\alpha$ and the number of edges in $S_i \cup \{v_ip(v_i)\}$ with label {\it blue} is 0 or 1 or $\beta$.
If the algorithm runs the {\it for loop} completely, then we are sure that the tree can be partitioned and
if at step $i$, there is no labeling for  $S_i$, then it shows that the label of edge $v_ip(v_i)$ should not be $f(v_ip(v_i))$.
So, the algorithm saves the correct label of $v_ip(v_i)$ in $g$, erases the labels of edges, breaks the {\it for loop} and starts the {\it for loop} from the beginning.
In the next iteration of the {\it for loop}, the algorithm has some restrictions in its labeling (if the label of an edge $e$ in $g(e)$ is not {\it free}, then its label must be equal to $g(e)$). Therefore, at step $i$ of the {\it for loop}, the algorithm should find a labeling like $h$ for the set of edges $S_i$ such that $|\{e: e\in S_i, h(e)=red\}\cup\{v_ip(v_i): f(v_ip(v_i))=red\}|\in \{0 , 1 ,\alpha\}$, also $|\{e: e\in S_i, h(e)=blue\}\cup\{v_ip(v_i): f(v_i)=blue\}|\in \{0 , 1 ,\beta\}$ and for each edge $e\in S_i$, if $g(e)\neq free$, then $g(e)=h(e)$. If the algorithm runs the {\it for loop} completely, then we are sure that the tree can be partitioned and
if at step $i$, there is no labeling for  $S_i$, then it shows that the label of the edge $v_ip(v_i)$ should not be $f(v_ip(v_i))$.
Now, if $g(v_ip(v_i))=f(v_ip(v_i))$, it shows that $T$ cannot be partitioned into two weakly semiregular forests and if $g(v_ip(v_i))=free$, the algorithm saves the correct label of $v_ip(v_i)$ in $g$, erases the labels, breaks the {\it for loop} and starts the {\it for loop} from the beginning. Note that if the algorithm does not run the {\it for loop} completely, then  the label of one edge in the function $g$ will be changed from $free$ into $red$ or $blue$. Thus, the {\it while loop} will be run at most $|E(T)|$ times. So, finally the algorithm finds a labeling or terminates and returns NO.

{\bf Complexity of Algorithm 1}\\
In Algorithm 1, if the {\it for loop} in Line 8, was completely run (if it was not {\it broken} in Line 16, 21 or 11), then
the algorithm will return YES in line 30. Otherwise, the label of one edge in function $g$ will be changed from $free$ into $red$ or $blue$. Thus, the {\it while loop} will be run at most $|E(T)|$ times. On the other hand, the  {\it for loop} in Line 8, takes  at most $ \mathcal{O}(n)$ times. Consequently, the running time of Algorithm 1 is $\mathcal{O}(n^2)$, hence we can determine whether $T$ can be partitioned into two weakly semiregular forests $T_1$ and $T_2$ in $\mathcal{O}(n^2)$.

(ii) Let $T$ be an arbitrary tree and $c$ be a constant number. Any subgraph of a tree is a forest, so if $T$ can be partitioned into $c$ weakly semiregular forests $T_1,T_2,  \ldots, T_c$, then there are $c$ numbers $\alpha_1,\alpha_2,\ldots,\alpha_c $ (not necessary distinct) such that $T_i$ is $(1,\alpha_i)$-forest. Without loss of generality, we can assume that $1 \leq \alpha_1\leq \alpha_2\leq \ldots\leq \alpha_c \leq \Delta(T) \leq n$. For each possible candidate for $(\alpha_1,\alpha_2,\ldots,\alpha_c )$
we run Algorithm 2. Note that the number of candidates is a polynomial in terms of the number of vertices (in terms of $n$).

Here, let us to introduce some notation and state a few properties of Algorithm 2.
Suppose that $|V(T)| =n$ and choose an arbitrary vertex $v\in V(T)$ to be its root.
Perform a breadth-first search
algorithm from the vertex $v$. This defines a partition $L_{0}, L_{1}, \ldots , L_{h}$ of the vertices of $T$
where each part $L_i$ contains the vertices of $T$ which are at depth $i$ (at distance exactly $i$
from $v$). Let $p(x)$ denote the neighbor of $x$ on the $xv$-path, i.e. its parent. Also,
 let $\lbrace v_{1}=v,v_2,\ldots , v_{n} \rbrace  $ be a list of the vertices according to their distance from the root.
If Algorithm 2 returns NO, it means that $T$ cannot be partitioned and if it returns YES, it means $T$ can be decomposed.
See Algorithm 2.

\begin{algorithm}
\algsetup{linenosize=\tiny}
\small
\caption{}
\begin{algorithmic}[1]
\STATE {{\bf Input:}} The tree $T$ and $c$ numbers $\alpha_1,\alpha_2,\ldots,\alpha_c $.
\STATE {{\bf Output:}}  Can $T$  be partitioned into
$c$ weakly semiregular forests $T_1,T_2,  \ldots, T_c$, such that $T_i$ is $(1,\alpha_i)$-forest.
\STATE Let $f:E(T) \rightarrow \{1,2,\ldots, c, free\}$ and put $f(e)\leftarrow free$ for all edges
\STATE Let $\ell:E(T) \rightarrow \{ 2^{\{1,2,\ldots,c\}}\}$ and put $\ell(e)\leftarrow \emptyset$ for all edges
\WHILE{ there is an edge $e$ such that $\ell(e)\neq \{1,2,\ldots,c\}$}
  \STATE For any edge $e$, put $f(e)\leftarrow free$
  \STATE $s\leftarrow $YES
  \FOR{ $i=1$ to $i=n$}
    \IF {there is no labeling like $h$ for the set of edges $S_i=\{v_iv_j: j>i\}$ with the colors $\{1,2,\ldots,c\}$ such that for each $k$, $1 \leq k \leq c$, $|\{e: e\in S_i, h(e)=k\}\cup\{v_ip(v_i):   f(v_ip(v_i))=k, i>1\}|\in \{0 , 1 ,\alpha_k\}$
    and for each edge $e\in S_i$ and color $t$, if $t\in \ell(e)$, then $h(e)\neq t$}
        \IF{ $f(v_ip(v_i))\in \ell(v_ip(v_i))$}
        \RETURN NO
        \ENDIF
         \IF{ $f(v_ip(v_i))\notin \ell(v_ip(v_i))$}
          \STATE $\ell(v_ip(v_i))\leftarrow \ell(v_ip(v_i)) \cup \{f(v_ip(v_i))\}$
           \STATE $s\leftarrow $NO
          \STATE {\bf break} the for loop
        \ENDIF
     \ENDIF
     \IF{$s=$YES}
    \STATE Label the set of edges $S_i=\{v_iv_j: j>i\}$ with the colors $\{1,2,\ldots,c\}$ such that for each $k$, $1 \leq k \leq c$, $|\{e: e\in S_i, h(e)=k\}\cup\{v_ip(v_i):   f(v_ip(v_i))=k, i>1 \}|\in \{0 , 1 ,\alpha_k\}$
    and for each edge $e\in S_i$ and color $t$, if $t\in \ell(e)$, then $h(e)\neq t$
    \STATE For each edge $e$ in $\{v_iv_j: j>i\}$ assign the label of $e$ to the variable $f(e)$
    \ENDIF
  \ENDFOR
  \IF{$s=$YES}
   \RETURN YES
  \ENDIF
\ENDWHILE
\RETURN NO
\end{algorithmic}
\end{algorithm}

In Algorithm 2, at each step the variable $\ell$ for each edge shows the set of forbidden colors for that edge. In other words, at any time, for each edge $e$, $\ell(e)$ is a subset of $\{1,2, \ldots,c\}$.

Assume that we want to find a labeling like $h$ with the labels $\{1,2,\ldots,c\}$ for the set of edges incident with a vertex $u$ such that for each $k$, $1 \leq k \leq c$, $|\{e: e\ni u, h(e)=k\}|\in \{\gamma_k\}$ and for each edge $e$ and color $t$, if $t\in \ell(e)$, then $h(e)\neq t$.
We claim that  this problem can be solved in polynomial time. Construct the bipartite graph $H$ with the vertex set $V(H)=\{a_1,a_2,\ldots,a_c\}\cup\{e: e\ni u\}$ such that $a_te\in E(H)$ if and only if $t\notin h(e)$. In this graph we want to find an $F$-factor such that
for each $i$, $F(a_i)=\gamma_k$ and for each edge $e$, $F(e)=1$. If the graph $H$ has an $F$-factor then there is a labeling like $h$ for the
the set of edges incident with the vertex $u$ with the specified properties.
In 1985, Anstee  gave a polynomial time algorithm for the $F$-factor problem and his algorithm
either returns one of the factors in question or shows that none exists, in $\mathcal{O}(n^3)$ time \cite{Anstee}.
Thus, Line 9 and Line 20 can be performed in polynomial time ({\bf Fact 1}).

{\bf Complexity of Algorithm 2}\\
In Algorithm 2, if the {\it for loop} in Line 8, was completely run (it was not {\it broken} in Line 16, 21 or 11), then
Algorithm will return YES in line 30. Otherwise, a color will be added to the set of forbidden colors of an edge. Thus, the {\it while loop} will be run at most $c|E(T)|$ times. On the other hand, by Fact 1, the  {\it for loop} in Line 8, takes  a polynomial  time to run. Consequently, the running time of Algorithm 1 is a polynomial in terms of $n$.

{\bf Sketch of Algorithm 2}\\
Performance of Algorithm 2 is similar to Algorithm 1, except that in Algorithm 2, at each step the variable $\ell$ for each edge shows the set of forbidden colors for that edge. So at any time, for each edge $e$, $\ell(e)$ is a subset of $\{1,2,3,\ldots,c\}$.
This completes the proof.
\\ \\
(iii) Suppose that $|V(T)| =n$ and choose an arbitrary vertex $v\in V(T)$ to be its root.
Perform a breadth-first search algorithm from the vertex $v$. This defines a partition $L_{0}, L_{1}, \ldots , L_{h}$ of the vertices of $T$
where each part $L_i$ contains the vertices of $T$ which are at depth $i$ (at distance exactly $i$
from $v$). Let $\lbrace v_{1}=v,v_2\ldots , v_{n} \rbrace  $ be a list of the vertices according to their distance from the root.
Do Algorithm 3 for the vertices of $T$ according to their indices.

\begin{algorithm}
\algsetup{linenosize=\tiny}
\small
\caption{}
\begin{algorithmic}[1]
\STATE {{\bf Input:}} The tree $T$.
\STATE {{\bf Output:}}  A decomposition of $T$ into $2\log_2 \Delta(T) + \mathcal{O}(1)$ weakly semiregular subgraphs.
\FOR{ $i=1$ to $i=n$}
     \STATE If $i=1$, let $a[0]a[1] \cdots a[\lfloor \log_2 \Delta \rfloor]$ be the binary representation of the number $d(v_1)$
(note that for each $r$, $a[r]$   is either a 1 or a 0)
and label the set of  edges $\{v_1v_j: j > 1\} $,
such that for each $t$, $0 \leq t \leq \lfloor \log_2 \Delta \rfloor$, if $a[t]=1$, then
 the number of edges incident with $v_1$
with label $t$ is $2^t$.

    \STATE If $i>1$, $v_i\in L_k$ and $k=0 (\mod 2)$. Let $a[0]a[1] \cdots a[\lfloor \log_2 \Delta \rfloor]$ be the binary representation of the number $d(v_i)-1$ and label the set of  edges $\{v_iv_j: j > i\} $,
such that for each $t$, $0 \leq t \leq \lfloor \log_2 \Delta \rfloor$, if $a[t]=1$, then
 the number of edges in $\{v_iv_j: j > i\}$
with label $t+ \lfloor \log_2 \Delta \rfloor+1$ is $2^t$.

 \STATE  If $i>1$, $v_i\in L_k$ and $k=1 (\mod 2)$. Let $a[0]a[1] \cdots a[\lfloor \log_2 \Delta \rfloor]$ be the binary representation of the number $d(v_i)-1$ and label the set of  edges $\{v_iv_j: j > i\} $,
such that for each $t$, $0 \leq t \leq \lfloor \log_2 \Delta \rfloor$, if $a[t]=1$, then
 the number of edges  in $\{v_iv_j: j > i\}$
with label $t $ is $2^t$.

\ENDFOR
\end{algorithmic}
\end{algorithm}

In Algorithm 3, for each $t$, $0 \leq t \leq \lfloor \log_2 \Delta \rfloor$, the set of edges with label $t$ forms a $(1,2^t)$-graph. Also, the set of edges with label $t+ \lfloor \log_2 \Delta \rfloor+1$ forms a $(1,2^t)$-graph.
Thus, one can see that the above-mentioned labeling is partitioning of edges into $2\log_2 \Delta(T) + \mathcal{O}(1)$ weakly semiregular subgraphs. This completes the proof.

\section{Proof of Theorem \ref{T2}}

It was shown in  \cite{dehghan2014complexity} that the following version of
 {\em Not-All-Equal (NAE) satisfying assignment problem} is {\bf NP}-complete.

 {\em Cubic Monotone  NAE (2,3)-Sat.}\\
\textsc{Instance}: Set $X$ of variables, collection $C$ of clauses over $X$ such that each
clause $c \in C$ has $\mid c  \mid \in \{ 2, 3\}$, every variable appears in
exactly three clauses and there is no negation in the formula.\\
\textsc{Question}: Is there a truth assignment for $X$ such that each clause in $C$ has at
least one true literal and at least one false literal?\\

We reduce {\em Cubic Monotone NAE (2,3)-Sat } to
our problem in polynomial time.
Consider an instance $ \Phi $, we transform this into a bipartite  graph $G$ in polynomial time
such that $ wr(G)=2$ if and only if $\Phi$ has an   NAE truth assignment.
We use three auxiliary gadgets $ \mathcal{H}_c$, $\mathcal{I}_c$ and $\mathcal{B}$ which are shown in Figure 1 and Figure 2.

\begin{figure}[ht]
\begin{center}
\includegraphics[scale=.45]{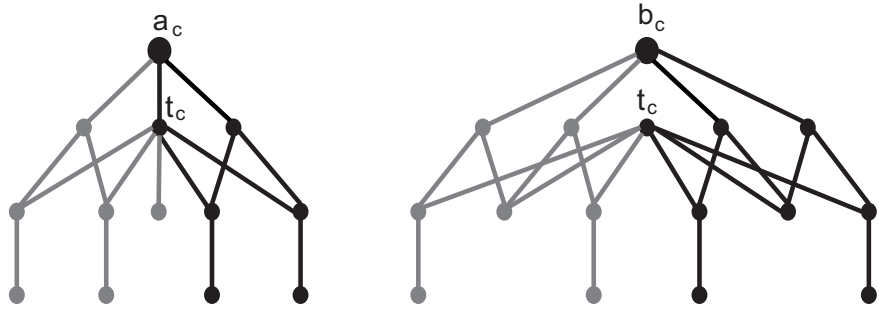}
\caption{The two gadgets $ \mathcal{H}_c$ and $\mathcal{I}_c$. The graph $\mathcal{I}_c$ is on the right hand side of
the figure.
}
\end{center}
\end{figure}

Our construction consists of three steps.\\
{\bf Step 1}. Put a copy of the graph $\mathcal{B}$, a copy of the complete bipartite graph $K_{1,6}$ and  a copy of the complete bipartite graph $K_{3,6}$.\\
{\bf Step 2}.  For each clause $c \in C$ with $\mid c \mid =3$, put
a copy of the gadget $\mathcal{H}_c $
and  for each clause $c \in C$ with $\mid c \mid =2$, put a copy
of the  gadget $\mathcal{I}_c $.\\
{\bf Step 3}.  For each variable $x \in X$, put a  vertex $x$ and for each clause $c =y \vee z \vee w$, where  $y,z,w \in X  $ add the edges $a_c y $, $a_c z  $ and $a_c w $. Also,
for each clause $c =y \vee z $, where  $y,z \in X  $ add the
edges $b_c y $ and $b_c z$.

Call the resultant graph $G$. The degree of every vertex in the graph $G$ is $ 1$ or $3$ or $6$ and the resultant graph is bipartite.
Assume that the graph $G$ can be partitioned into two weakly semiregular graphs $G_1$ and $G_2$, we have the following lemmas.
\\ \\
{\bf Lemma 1.} The graphs $G_1$ and  $G_2$ are (1,3)-graph.
\\ \\
{\bf Proof of Lemma 1.} Without loss of generality, assume that the graph $G_1$ is $ (\alpha_1,\alpha_2)$-graph and
the graph $G_2$ is $ (\beta_1,\beta_2)$-graph.
Since $\Delta(G)=6$, by attention to the structure of the graph $\mathcal{B}$, with respect to the symmetry, the following cases for $(\alpha_1\alpha_2,\beta_1\beta_2)$ can be considered: $(16,12),(15,12),(24,12),(14,12),(13,13)$.
The graph $G$ contains a copy of the complete bipartite graph $K_{1,6}$, so the case (24,12) is not possible, also,  the graph $G$ contains  a copy of complete bipartite graph $K_{3,6}$, so the cases $(16,12),(15,12),(14,12)$ are not possible.
Thus, the graphs $G_1$ and  $G_2$ are (1,3)-graph. $ \diamondsuit$
\\ \\
{\bf Lemma 2.} For every vertex $v$ with degree three all edges incident with the vertex
 $v$ are in one part.
\\ \\
{\bf Proof of Lemma 2.} Since the graphs $G_1$ and  $G_2$ are (1,3)-graph, the proof is clear. $ \diamondsuit$

Now, we present an NAE truth assignment for the formula $\Phi$.
For every $x\in X$, if all edges incident with the vertex $x$ are in the graph $G_1$, put $\Gamma(x)=true$ and if all edges incident with
the vertex $x$ are in graph $G_2$,
put $\Gamma(x)=false$.
Let $c =y \vee z \vee w$ be an arbitrary clause, if all edges $a_c y ,a_c z ,a_c w $ are in the graph $G_1$ ($G_2$, respectively),
then by the construction of the gadget $\mathcal{H}_c$, the degree of the  vertex $t_c $ in the graph $G_2$ ($G_1$, respectively) is at least 5 (5, respectively). This is a contradiction. Similarly, let $c =y \vee z  $ be an arbitrary clause, if all edges $b_c y ,b_c z  $ are in the graph $G_1$ ($G_2$, respectively),
then by the construction of $\mathcal{I}_c$, the degree of the  vertex $t_c $ in the graph $G_2$ ($G_1$, respectively) is 6. This is a contradiction.
Hence, $\Gamma$ is an NAE satisfying assignment.
On the other hand, suppose that $\Phi$ is  NAE satisfiable with the satisfying
assignment $\Gamma : X  \rightarrow \{true, false\}$. For every variable $x\in X$, put all
edges incident with the the vertex $x$ in $G_1$ if and only if $\Gamma (x)=true$. By this method,
one can show that the graph $G$ can be partitioned into two weakly semiregular subgraphs.
This completes the proof.

\begin{figure}[ht]
\begin{center}
\includegraphics[scale=.30]{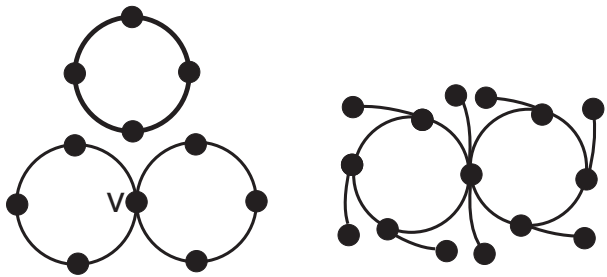}
\caption{The two gadgets $\mathcal{P}$ and $ \mathcal{B}$. The graph $\mathcal{B}$ is on the right hand side of
the figure.
}
\end{center}
\end{figure}

\section{Proof of Theorem \ref{T3}}

(i) Any subgraph of a tree is a forest, so in every decomposition of a tree $T$ into some semiregular subgraphs, each subgraph is a (1,2)-forest. Thus, $sr(T) \geq \lceil\frac{\Delta(T)}{2}\rceil$. For every bipartite  graph $H$, we have $\chi'(H)=\Delta(H)$ (see for example \cite{MR1367739}). Assume that $f:E(T) \rightarrow \{1,\ldots, \Delta(T)\}$ is a proper edge coloring for $T$. The following partition is a decomposition of $T$ into $\lceil\frac{\Delta(T)}{2}\rceil$  semiregular subgraphs.

$$E(T)=\bigcup_{i=1}^{\lceil\frac{\Delta(T)}{2}\rceil}\{e: f(e)=i \text{ or } f(e)=i +\lceil\frac{\Delta(T)}{2}\rceil\}. $$
\\
This completes the proof.

(ii)  Let $G$ be an arbitrary graph. By Vizing's theorem \cite{MR0180505}, the edge chromatic number of a graph $G$ is equal to either $ \Delta(G) $ or $ \Delta(G) +1 $. The following partition is a decomposition of $G$ into $\lceil\frac{\chi'(G)}{2}\rceil$  semiregular subgraphs.

$$E(G)=\bigcup_{i=1}^{\lceil\frac{\chi'(G)}{2}\rceil}\{e: f(e)=i \text{ or } f(e)=i +\lceil\frac{\chi'(G)}{2}\rceil\}. $$
\\
So the graph $G$ can be partitioned into $\lceil\frac{\Delta(G)+1}{2}\rceil$ semiregular subgraphs.

(iii) We use a reduction from the following {\bf NP}-complete problem \cite{dehghan2014complexity}.

 {\em Cubic Monotone  NAE (2,3)-Sat.}\\
\textsc{Instance}: Set $X$ of variables, collection $C$ of clauses over $X$ such that each
clause $c \in C$ has $\mid c  \mid \in \{ 2, 3\}$, every variable appears in
exactly three clauses and there is no negation in the formula.\\
\textsc{Question}: Is there a truth assignment for $X$ such that each clause in $C$ has at
least one true literal and at least one false literal?\\

We reduce {\em Cubic Monotone NAE (2,3)-Sat } to
our problem in polynomial time.
Consider an instance $ \Phi $, we transform this into a bipartite  graph $G$ with $\Delta(G)\leq 6$ in polynomial time
such that $ sr(G)=2$ if and only if $\Phi$ has an   NAE truth assignment.
We use three auxiliary gadgets $ \mathcal{D}_c$,  $\mathcal{F}_c$  and $\mathcal{P}$ which are shown in Figure 3 and Figure 2.

\begin{figure}[ht]
\begin{center}
\includegraphics[scale=.45]{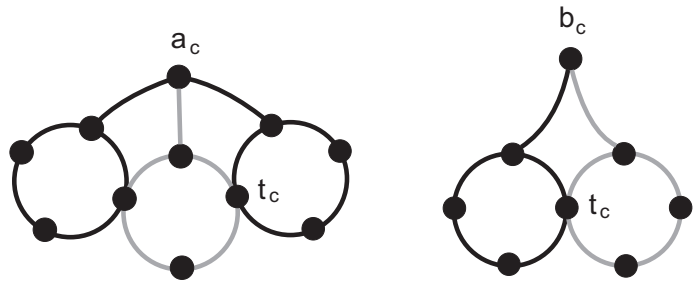}
\caption{The two auxiliary gadgets $ \mathcal{F}_c$ and $\mathcal{D}_c$. $\mathcal{D}_c$ is on the right hand side of
the figure.
}
\end{center}
\end{figure}

Our construction consists of three steps.\\
{\bf Step 1}. Put a copy of the graph $\mathcal{P}$.\\
{\bf Step 2}. For each clause $c \in C$ with $\mid c \mid =3$, put
a copy of the gadget $\mathcal{F}_c $
and  for each clause $c \in C$ with $\mid c \mid =2$, put a copy
of the gadget $\mathcal{D}_c $.\\
{\bf Step 3}.  For each variable $x \in X$, put a  vertex $x$ and for each clause $c =y \vee z \vee w$, where  $y,z,w \in X  $ add the edges $a_c y $, $a_c z  $ and $a_c w $. Also,
for each clause $c =y \vee z $, where  $y,z \in X  $ add the
edges $b_c y $ and $b_c z$.

Call the resultant graph $G$. The degree set of the graph $G$ is $\{ 2,3,4,6\}$ and the graph is bipartite.
Assume that $G$ can be partitioned into two semiregular graphs $G_1$ and $G_2$, we have the following lemmas.
\\ \\
{\bf Lemma 3.} The graphs $G_1$ and  $G_2$ are (2,3)-graph.
\\ \\
{\bf Proof of Lemma 3.} Without loss of generality assume that $G_1$ is $(\alpha-1,\alpha)$-graph such that $\Delta(G_1)=\alpha$ and $G_2$ is $(\beta-1,\beta)$-graph such that $\Delta(G_2)=\beta$.
In the graph $G$ any vertex of degree six has a neighbor of degree three, Thus, $\alpha \neq 6$ and $\beta\neq 6$.
Also, there is no vertex of degree five, and any neighbor of each vertex of degree six has degree three, so we can assume that  $\alpha \neq 5$ and $\beta\neq 5$.
In the graph $\mathcal{P}$, the degree of the vertex $v$ is four and the degree of each of its neighbor is two (note that the graph $G$ contains a copy of the graph $\mathcal{P}$). Thus, by the structure of $\mathcal{P}$ and since the graph $G$ contains a vertex of degree six, we have $\alpha \neq 4$ and $\beta\neq 4$.
On the other hand, since $\Delta(G)=6$, we have $\alpha=\beta=3$. Hence, the graphs $G_1$ and  $G_2$ are (2,3)-graph. $ \diamondsuit$
\\ \\
{\bf Lemma 4.} For every vertex $z$ with degree three or two all edges incident with the vertex
 $z$ are in one part.
\\ \\
{\bf Proof of Lemma 4.} Since the graphs $G_1$ and  $G_2$ are (2,3)-graph, the proof is clear. $ \diamondsuit$

Now, we present an NAE truth assignment for the formula $\Phi$.
For every $x\in X$, if all edges incident with the vertex $x$ are in $G_1$, put $\Gamma(x)=true$ and if all edges incident with
the vertex $x$ are in $G_2$,
put $\Gamma(x)=false$.
Let $c =y \vee z \vee w$ be an arbitrary clause, if all edges $a_c y ,a_c z ,a_c w $ are in $G_1$ ($G_2$, respectively),
then by the construction of the  gadget $\mathcal{F}_c$ and Lemma 4, the degree of the vertex $t_c $ in the graph $G_2$ ($G_1$, respectively) is 4 (4, respectively). This is a contradiction. Similarly, let $c =y \vee z  $ be an arbitrary clause, if all edges $b_c y ,b_c z  $ are in $G_1$ ($G_2$, respectively),
then by the construction of the gadget $\mathcal{D}_c$, the degree of the vertex $t_c $ in the graph $G_2$ ($G_1$, respectively) is 4. This is a contradiction.
Hence, $\Gamma$ is an NAE satisfying assignment.
On the other hand, suppose that $\Phi$ is  NAE satisfiable with the satisfying
assignment $\Gamma : X  \rightarrow \{true, false\}$. For every variable $x\in X$, put all
edges incident with the the vertex $x$ in $G_1$ if and only if $\Gamma (x)=true$. By this method,
it is easy to show that $G$ can be partitioned into two  semiregular subgraphs.
This completes the proof.

\section{Proof of Theorem \ref{T4}}

Let $G$ be a graph.
We say that an edge-labeling $\ell:E(G)\rightarrow \mathbb{N}$ is an additive vertex-colorings if and only if for each edge
$uv$, the sum of labels of the edges incident to $u$ is different
from the sum of labels of the edges incident to $v$.
It was shown that
 determining whether a given 3-regular graph $G$ has
an edge-labeling which
is an additive vertex-coloring from $\{1,2\}$  is
$ \mathbf{NP} $-complete \cite{MR3072733}. For a given 3-regular graph $G$, it is easy to see that $G$  has an edge-labeling which
is an additive vertex-coloring from $\{1,2\}$ if and only if the edge set of $G$ can be partitioned
 into at most two locally irregular subgraphs. Thus,  determining whether a given 3-regular graph $G$ can be decomposed into two locally irregular subgraphs is {\bf NP}-complete \cite{MR3072733}. We will reduce this problem to our problem.
Let $G$ be a 3-regular graph. We construct a graph $G'$ such that $G$ can be partitioned into two locally irregular subgraphs if and only if $G'$
can be partitioned into 2 subgraphs, such that each subgraph
is locally irregular or weakly semiregular. Let $G'= G \cup C_4 \cup P_5 \cup K_{9,9} \bigcup_{i=4}^{8} K_{1,i} $.

The degree set of $G'$ is $D=\{j: 1\leq j \leq 9\}$, so $|D|\geq 9$, thus $G'$ cannot be partitioned into two weakly semiregular subgraphs. Now, assume that $G'$ can be partitioned into two subgraphs $\mathcal{I}$ and $\mathcal{R}$ such that $\mathcal{I}$ is locally irregular and $\mathcal{R}$ is $(\alpha, \beta)$-graph. The graph $G'$ contains a copy of $C_4$, thus $2\in \{\alpha,\beta\}$. Also, $G'$ contains a copy of $P_5$, thus $1\in \{\alpha,\beta\}$. Hence $\mathcal{R}$ is  a $(1,2)$-graph. Note that $ K_{9,9}$ cannot be partitioned into two subgraphs $\mathcal{I}$ and $\mathcal{R}$ such that $\mathcal{I}$ is locally irregular and $\mathcal{R}$ is $(1,2)$-graph. Thus, $G'$ cannot be partitioned into two subgraphs $\mathcal{I}$ and $\mathcal{R}$ such that $\mathcal{I}$ is locally irregular and $\mathcal{R}$ is weakly semiregular.
On the other hand, it is easy to see that the graph $ C_4 \cup P_5 \cup K_{9,9} \bigcup_{i=4}^{8} K_{1,i}$ can be partitioned into two locally irregular subgraphs. Therefore, $G$ can be partitioned into two locally irregular subgraphs if and only if $G'$
can be partitioned into 2 subgraphs, such that each subgraph
is locally irregular or weakly semiregular. This completes the proof.

\section{Proof of Theorem \ref{T5}}
(i) Let $\epsilon>0$ be a fixed number and
$G$ be a 3-regular graph with sufficiently large number of
vertices in terms of $\epsilon$.
Construct the graph $H$ from the graph $G$ by replacing every edge $ab$ of $G$ by a copy of the gadget $I(a,b)$ which is shown in Fig. 4.
It was shown that it is $\mathbf{NP}$-complete to determine whether the edge chromatic number of a
cubic graph is three \cite{MR635430}. Assume that the number of vertices in the graph $H$ is $n$.
We show that if
$\chi'(G)=3$ then $rep(H^{c})\leq (1+\epsilon)(\frac{n}{2})^{3}$
and if $\chi'(G)>3$ then $rep(H^{c})\geq(\frac{n}{2})^{4}$,
consequently, there is no polynomial time $\theta$-approximation algorithm for
computing $rep(A^{c})$ when
$$\frac{(\frac{n}{2})^{4}}{(1+\epsilon)(\frac{n}{2})^{3}}=\frac{n}{2(1+\epsilon)}>
(1-\epsilon)\frac{n}{2}=\theta.$$ By the structure of the gadget $I(a,b)$, the graph $H$ is 3-regular and
triangle-free, also by the structure of $H$, $\chi'(G)=\chi'(H)$. Note that if $H$ is triangle-free then $rep(H^{c})$ is
square free. Let $rep(H^{c})=\prod_{i=1}^{d}p_{i}$, where for each $i$, $i=1,\ldots,d$, $p_{i}$ is a prime number. Assume that
$r:V(H^{c})\longrightarrow \mathbb{Z}_{rep(H^{c})}$ is a
representation for $H^{c}$. For each vertex $v$ of $H^{c}$ define a $d$-triple $(r_{v}^{1},\ldots,r_{v}^{d})\in
\prod_{i=1}^{d}\mathbb{Z}_{p_{i}}$
such that $r_{v}^{i}= (r(v) \mod p_{i})$. By the definition of the function $r$,  $vw$ is
an edge in $H$ if and only if there exists an index $i$ such that $r^{i}_{v}=r^{i}_{w}$. For
each edge $e=xy$ of $H$, define $S(e)=\{i:r^{i}_{x}=r^{i}_{y}\}$.
So for each edge $e$, $S(e)$ is non-empty and since the graph $H$ is triangle-free, for every two incident edges $e$ and $e'$ we
have $S(e)\cap S(e')=\emptyset$. Let $c:E(H)\longrightarrow
\{1,\ldots,d\}$ be a function such that $c(e)=\min S(e)$. It is clear that
$c$ is a proper edge coloring for the graph $H$. So
\begin{equation} \label{1}
d\geq \chi'(H)
\end{equation}

\begin{figure}[ht]
\begin{center}
\includegraphics[scale=.3]{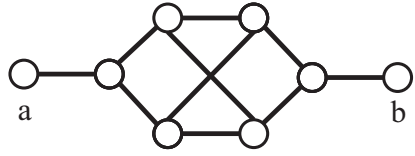}
\caption{The auxiliary graph $I(a,b)$.
}
\end{center}
\end{figure}

Define $M_{i}=\{e\in E(H), i\in S(e)\}$ for every $i$, $1 \leq i \leq d$.
The set $M_{i}$ contains all edges of the graph $H$ like $e=vu$ such that $v$ and $u$ have a same $i$-th
component. Since the graph $H$ is triangle-free, it follows that the set of edges $M_{i}$ is
a matching. Also, each $z\in \mathbb{Z}_{p_{i}}$ appears at most 2
times as the $i$-th component of vertices in the graph $H$. Also, each vertex of
$H$ which is not adjacent to any vertex of $M_{i}$, has a unique
$i$-th component. For each $i$ denote the number of
edges of $M_{i}$ by   $m_{i}$ (note that $|V(H)|=|V(H^c)|=n$). We have:
\begin{equation} \label{2}
p_{i}\geq m_{i}+(n-2m_{i})=n-m_{i}
\end{equation}
Also, since every matching
has at most $\frac{n}{2}$ edges, it follows that
\begin{equation} \label{3}
p_{i}\geq \frac{n}{2}
\end{equation}
Now, let $\chi'(H)=3$ and $f:E(H) \longrightarrow \{1,2,3\}$ be a
proper edge coloring of $H$. The edges of $H$ can be partitioned into three
perfect matchings $f_1,f_2$ and $f_3$, where $f_i=\{e:f(e)=i\}$.
Without loss of generality for each $i$, $i=1,2,3$, assume that
$f_{i}=\{e_{i}^{1},\ldots,e_{i}^{\frac{n}{2}}\}$.

It follows from the prime number theorem that for any real $\alpha >0$ there is a $n_0 > 0 $ such that for all $n' > n_0$ there is a prime $p$ such that $n' < p < (1 + \alpha) n'$ (see \cite{PNT}  page 494).
Thus for a sufficiently large number $n$, there are three prime
numbers $p_1,p_2,p_3$ such that $\frac{n}{2} \leq  p_{1}< p_{2} < p_{3}<
\frac{n}{2}(1+\epsilon)^{\frac{1}{3}} $. For every vertex $v$ of the graph $H$,
call the  edges  incident with the vertex $v$,  $e_{1}^{\alpha}$,
$e_{2}^{\beta}$ and $e_{3}^{\gamma}$ and let
$\psi:V(H^{c})\longrightarrow \mathbb{Z}_{p_{1}p_{2}p_{3}}$ be
a function such that  $\psi(v)=(\alpha,\beta,\gamma)$. Clearly, this is a
representation, so $rep(H^{c})\leq p_{1}p_{2}p_{3} <
(1+\epsilon)(\frac{n}{2})^{3}$. On the other side, assume that
$\chi'(G)>3$, so $\chi'(H)>3$. Thus, we have:

\begin{align*}
rep(H^{c}) & =     \prod_{i=1}^{d} p_{i}                       \\
           & \geq  \prod_{i=1}^{4} p_{i}                        \,\,& {\text{By Equation \ref{1},}}\\
           & \geq      (\frac{n}{2})^{4}                         \,\,& {\text{By Equation \ref{3},}}
\end{align*}
This completes the proof.
\\
 \\
(ii)
Let $G^c$ be a triangle-free $r$-regular
graph.  By Vizing's theorem \cite{MR0180505}, the chromatic index of $G^c$ is equal to either $ \Delta(G^c) $ or $ \Delta(G^c) +1 $. Thus, for every $r$-regular graph $G^c$, $r \leq \chi'(G^c) \leq r+1$.
Therefore, the set of edges of $G^c$ can be partitioned into $r+1$ matchings
 $M_{1},\ldots, M_{r+1}$.  By an argument similar to the proof of part (i),  we have:

\begin{align*}
rep(G) & \leq (1+\epsilon)\prod_{i=1}^{r+1}(n-|M_{i}|) \\
       & \leq (1+\epsilon)(n-\frac{rn}{2(r+1)})^{r+1}\,\,& {\text{By Equation \ref{2},}}\\
       & \leq (1+\epsilon)e(\frac{n}{2})^{r+1}        \,\,\,\,&{\text{By inequality }(1+\dfrac{1}{x})^x <e,}\\
        \text{On the other hand}\\
rep(G) & \geq (\frac{n}{2})^{r}.
\end{align*}

Therefore we have a polynomial time $(1+\epsilon)\frac{e}{2}n$
approximation algorithm for
computing $rep(G)$.

\section{Concluding remarks and future work}

\subsection{Trees}

We proved that for every tree $T$, $wr(T)\leq 2\log_2 \Delta(T) + \mathcal{O}(1)$.
On the other hand, there are infinitely many values of $\Delta$ for which the tree $T$ might be chosen so that $wr(T)\geq \log_3 \Delta(T)$. Finding the best upper bound for trees can be interesting.

\begin{prob}
Find the best upper bound for the weakly semiregular numbers of trees in terms of the maximum degree.
\end{prob}

We proved that there is an $\mathcal{O}(n^2)$ time algorithm to determine whether the weakly semiregular number of a given tree is two. Also, if $c$ is a constant number, then there is a polynomial time algorithm to determine whether the weakly semiregular number of a given tree is at most $c$.
However, one further step does not seem trivial. Is there any polynomial time algorithm to determine  the weakly semiregular number of trees?

\begin{prob}
Is there any polynomial time algorithm to determine  the weakly semiregular number of trees?
\end{prob}

In this work we present an algorithm with running time $\mathcal{O}(n^2)$ to determine whether the weakly semiregular number of a given tree is at most two. Is there any algorithm with running time $\mathcal{O}(n\lg n)$ for this problem?

\begin{prob}
Is there any algorithm with running time $o(n^2)$ for determining  whether the weakly semiregular number of a given tree is at most two?
\end{prob}

\subsection{Planar graphs}

Balogh {\it et al.}  proved that a planar graph can be partitioned into three forests so
that one of the forests has maximum degree at most 8 \cite{balogh2005covering}. On the other hand, we proved that
for every tree $T$, $wr(T)\leq 2\log_2 \Delta(T) + \mathcal{O}(1)$. Thus, for every planar graph $G$, we have
$wr(G)\leq 4\log_2 \Delta(G) + \mathcal{O}(1)$. Also,  it was shown that every planar graph with girth $g\geq 6$ has an edge partition into two forests, one having maximum degree at most 4 \cite{gonccalves2009covering}. Thus, for every planar graph
$G$ with girth $g\geq 6$, we have $wr(G)\leq 2\log_2 \Delta(G) + \mathcal{O}(1)$. Finding a good upper bound for all planar graphs can be interesting.

\begin{prob}
Is this true "For every planar graph $G$, we have
$wr(G)\leq 2\log_2 \Delta(G) + \mathcal{O}(1)$"?
\end{prob}

\subsection{Representation Number}

In this work, we proved that if $\mathbf{NP\neq P}$, then for any $\epsilon >0$, there is no polynomial time
$(1-\epsilon)\frac{n}{2}$-approximation algorithm for the computation of
representation number of regular graphs with $n$ vertices. In 2000 it was shown by Evans, Isaak and Narayan \cite{MR1782042} that if $n,m\geq 2$, then $rep(nK_m)=p_ip_{i+1}\ldots p_{i+m-1}$ where $p_i$ is the smallest prime number greater than or equal to $m$ if and only if there exists a set of $n-1$ mutually orthogonal Latin squares of order $m$.
It is interesting to investigate what
our result implies about the Orthogonal Latin Square Conjecture (There exists $n-1$ mutually orthogonal Latin
squares of order $n$ if and only if $n$ is a prime power). That is, can our reduction be extended from regular graphs
to just $nK_m$.

\section{Acknowledgments}
\label{}

The authors would like to thank the anonymous referees for their useful comments
 which helped to improve the presentation of this paper.

\bibliographystyle{plain}
\bibliography{luckyref}

\end{document}